%% file: P1P2.tex
\documentclass[11pt]{article}
\input{format} \input{mathdefs} \input{thmdefs}

\newcommand{\G}{\mathbb{G}}
\newcommand{\llangle}{\langle\langle}
\newcommand{\rrangle}{\rangle\rangle}
\newcommand{\Ball}{\mathrm{Ball}}
\title{Reiter's properties $(P_1)$ and $(P_2)$ \\ for locally compact quantum groups}
\author{\textit{Matthew Daws} \and \textit{Volker Runde}\thanks{Research supported by NSERC}}
\date{}
\begin{document}
\maketitle
\begin{abstract}
A locally compact group $G$ is amenable if and only if it has Reiter's property $(P_p)$ for $p=1$ or, equivalently, all $p \in [1,\infty)$, i.e., there is a net $( m_\alpha )_\alpha$ of non-negative norm one functions in $L^p(G)$ such that $\lim_\alpha \sup_{x \in K} \| L_{x^{-1}} m_\alpha - m_\alpha \|_p = 0$ for each compact subset $K \subset G$ ($L_{x^{-1}} m_\alpha$ stands for the left translate of $m_\alpha$ by $x^{-1}$). We extend the definitions of properties $(P_1)$ and $(P_2)$ from locally compact groups to locally compact quantum groups in the sense of J.\ Kustermans and S.\ Vaes. We show that a locally compact quantum group has $(P_1)$ if and only if it is amenable and that it has $(P_2)$ if and only if its dual quantum group is co-amenable. As a consequence, $(P_2)$ implies $(P_1)$.
\end{abstract}
\begin{keywords}
amenability, co-amenability, Leptin's theorem, locally compact quantum groups, operator spaces, Reiter's property $(P_1)$, Reiter's property $(P_2)$.
\end{keywords}
\begin{classification}
Primary 43A07; Secondary 22D35, 43A30, 46L07, 46L65, 46L89, 47L25, 47L50, 81R15, 81R50.
\end{classification}
\section*{Introduction}
A locally compact group $G$ is said to be \emph{amenable} if there is an invariant mean on $L^\infty(G)$, i.e., a state $M$ of the von Neumann algebra $L^\infty(G)$ such that
\[
  \langle L_x \phi, M \rangle = \langle \phi, M \rangle  \qquad (\phi \in L^\infty(G), \, x \in G).
\]
(If $f$ is any function on $G$ and $x \in G$, we denote by $L_x f$ the left translate of $f$ by $x$, i.e., $(L_x f)(y) := f(xy)$ for $y \in G$.) Approximating $M$ in the weak$^\ast$ topology of $L^\infty(G)^\ast$ by normal states, i.e., non-negative, norm one functions in $L^1(G)$ and then passing to convex combinations, we obtain a net $(m_\alpha )_\alpha$ of such functions in $L^1(G)$ that is asymptotically invariant in the sense that
\begin{equation} \label{Reiter}
  \lim_\alpha \| L_{x^{-1}} m_\alpha - m_\alpha \|_1 = 0 \qquad (x \in G).
\end{equation}
On the other hand, whenever we have a net $( m_\alpha )_\alpha$ of non-negative norm one functions in $L^1(G)$ satisfying (\ref{Reiter}), then each of its weak$^\ast$ accumulation points in $L^\infty(G)^\ast$ is a left invariant mean, so that $G$ is amenable.
\par 
Even though it is not obvious, the net $( m_\alpha )_\alpha$ can be chosen for amenable $G$ in such a way that the convergence in (\ref{Reiter}) is uniform in $x$ on each compact subset of $G$ (\cite[Proposition 6.12]{Pie}), a condition called \emph{Reiter's property $(P_1)$} in the literature. More generally, one can define \emph{Reiter's property $(P_p)$} for any $p \in [1,\infty)$ (\cite[Definition 8.3.1]{RSt}), but as it turns out, the properties $(P_p)$ are all equivalent (\cite[Theorem 8.3.2]{RSt}). In \cite{LoA}---see Secton \ref{prelude} below---, the equivalence of amenability, $(P_1)$, and $(P_2)$ was used to prove Leptin's theorem (\cite{Lep}): $G$ is amenable if and only if $A(G)$, Eymard's Fourier algebra (\cite{Eym}), has a bounded approximate identity.
\par 
Leptin's theorem assumes a very natural form in the language of Kac algebras (see \cite{ES2}). In this language---using the terminology of \cite{BT}---, Leptin's theorem reads as: a locally compact group $G$, if viewed as a Kac algebra, is amenable if and only if its Kac algebraic dual is co-amenable. Hence, it is only natural to ask whether Leptin's theorem holds true for arbitrary Kac algebras: a Kac algebra is amenable if and only if its dual is co-amenable. In \cite{Voi}, D.\ V.\ Voiculescu showed that, indeed, the co-amenability of a Kac algebra implies the amenability of its dual. In \cite{ES1}, it was claimed that the converse is also true, but the proof given in \cite{ES1} contains an error. Ultimately, Z.-J.\ Ruan was able to salvage the result at least for discrete Kac algebras (\cite{Rua2}) whereas the general case remains open.
\par
Recently, J.\ Kustermans and S.\ Vaes introduced a sur\-pri\-sing\-ly simple system of axioms for what they call \emph{locally compact quantum groups} (\cite{KV1} and \cite{KV2}): those axioms cover the Kac algebras (and therefore all locally compact groups), allow for the development of a Pontryagin type duality theory, but also seem to cover all known examples of $\cstar$-algebraic quantum groups, such as Woronowicz's $\mathrm{SU}_q(2)$ (\cite{Wor}). For a detailed exposition on the history of locally compact quantum groups---with many references to the original literature---, we refer to the introduction of \cite{KV1} and to \cite{Vai}. Of course, the question whether amenability is dual to co-amenability--- so that Leptin's theorem holds true for locally compact quantum groups--- is a natural one, and---as for Kac algebras---it is only known to be true in the discrete case (\cite{Tom}).
\par 
The problem to prove Leptin's theorem for general locally compact quantum groups appears to be formidable. R.\ Tomatsu's proof in the discrete case (see \cite{Tom}) makes heavy use of the particular structure of discrete quantum groups (as does Ruan's argument in the discrete Kac algebra case) and does not appear to be adaptable to the general locally compact situation.
\par 
The present paper grew out the attempt to extend the proof of Leptin's theorem from \cite{LoA} to locally compact quantum groups. The problems arising with such an endeavor are numerous. How can Reiter's properties $(P_1)$ and $(P_2)$ be formulated? How do $(P_1)$ and $(P_2)$ relate to amenability and co-amenability, respectively? Finally, are $(P_1)$ and $(P_2)$ equivalent?
\par 
We proceed as follows. The first two sections are mostly expository. We recall the definition of Reiter's properties and reformulate them in a way that will later allow us to extend them to a quantum group setting. Then we give a brief overview of locally compact quantum groups (with references to the original literature). With these preparations, we then define property $(P_1)$ for quantum groups and show that $(P_1)$ and amenability are indeed equivalent; both the definition of $(P_1)$ and the proof of the equivalence result rely heavily on the theory of operator spaces (\cite{ER}, \cite{Pis}, and \cite{Pau}). We then go on and define $(P_2)$ for quantum groups, and we show that $(P_2)$ is equivalent, not just to the amenability of the quantum group, but to the co-amenability of its dual (again, both the definition and the result are steeped in operator space theory). As a consequence, $(P_2)$ implies $(P_1)$ whereas the converse remains open.
\section{Leptin's theorem through $(P_1)$ and $(P_2)$} \label{prelude}
The original proof of Leptin's theorem, as given in \cite{Lep}, relied on F{\o}lner type conditions, for which it is difficult to see how---if at all---they can be transferred to the context of general locally compact quantum groups. In \cite{LoA}, an alternative proof---making use of properties $(P_1)$ and $(P_2)$ instead---was attempted, but the argument given in \cite{LoA} was incomplete. 
\par
We begin this section with recalling Reiter's properties $(P_p)$ for $p \in [1,\infty)$ (\cite[Definition 8.3.1]{RSt}):
\begin{definition} \label{Ppdef}
Let $G$ be a locally compact group, and let $p \in [1,\infty)$. We say that $G$ has \emph{Reiter's property $(P_p)$} if there is a net $( m_\alpha )_\alpha$ of non-negative norm one functions in $L^p(G)$ such that
\[
  \lim_\alpha \sup_{x \in K} \| L_{x^{-1}} m_\alpha - m_\alpha \|_p = 0
\]
for all compact $K \subset G$.
\end{definition}
\begin{remarks}
\item It is not difficult to see that $G$ has $(P_p)$ for all $p \in [1,\infty)$ if and only if it has $(P_1)$ (\cite[Theorem 8.3.2]{RSt}).
\item By \cite[Proposition 6.12]{Pie}, $(P_1)$ is equivalent to amenability.
\end{remarks}
\par 
We now indicate how the argument in \cite{LoA} can be repaired:
\begin{proof}[Proof of Leptin's theorem via properties $(P_1)$ and $(P_2)$]
Let $G$ be a locally compact group, and suppose that $G$ is amenable, i.e., has Reiter's property $(P_1)$ and thus, equivalently, $(P_2)$. This means that is a net $( \xi_\alpha )_{\alpha \in \mathbb A}$ of non-negative norm one functions in $L^2(G)$ such that 
\[
  \lim_\alpha \sup_{x \in K} \| L_{x^{-1}} \xi_\alpha - \xi_\alpha \|_2 = 0
\]
for all compact sets $K \subset G$. For $\alpha \in \mathbb A$, define
\[
  e_\alpha \!: G \to \comps, \quad x \mapsto \langle L_{x^{-1}} \xi_\alpha, \xi_\alpha \rangle.
\]
Then $( e_\alpha )_\alpha$ is a net in $A(G)$ converging to $1$ uniformly on all compact subsets of $G$. By \cite[Theorem B$_2$]{GL}, this is enough for $( e_\alpha )_\alpha$ to be a bounded approximate identity for $A(G)$.
\par 
The converse implication of Leptin's theorem is easier to prove (and has long been known to extend to locally compact quantum groups; see \cite{BT}).
\end{proof}
\par 
We conclude this section with a recasting of Definition \ref{Ppdef} that will enable us later to extend it from locally compact groups to quantum groups (at least for $p =1, 2$).
\par 
Let $G$ be a locally compact group, let $p \in [1,\infty)$, and let $g \in L^p(G)$. Then
\[
  L_\bullet(g) \!: G \to L^p(G), \quad x \mapsto L_{x^{-1}} g
\]
is a bounded, continuous function with values in $L^p(G)$. Let $f \in {\cal C}_0(G)$, and define $fL_\bullet(g) \!: G \to L^p(G)$ pointwise, i.e., $(fL_\bullet(g))(x) := f(x) L_{x^{-1}} g$ for $x \in G$. Since $f \in {\cal C}_0(G)$, $fL_\bullet(g)$ also vanishes at infinity and thus lies in ${\cal C}_0(G,L^p(G)) \cong {\cal C}_0(G) \tensor^\lambda L^p(G)$ (following \cite{ER}, we denote the injective Banach space tensor product by $\tensor^\lambda$). 
\par 
Let $( m_\alpha )_\alpha$ be a bounded net in $L^p(G)$. Then it is straightforward to verify that
\[
  \lim_\alpha \sup_{x \in K} \| L_{x^{-1}} m_\alpha - m_\alpha \|_p = 0
\]
holds for each compact $K \subset G$ if and only if
\[
  \lim_\alpha \| f L_\bullet(m_\alpha) - f \tensor m_\alpha \|_{{\cal C}_0(G) \tensor^\lambda L^p(G)} = 0
\]
is true for all $f \in {\cal C}_0(G)$.
\par 
In view of this, Definition \ref{Ppdef} and the following are equivalent:
\begin{definition} \label{Ppdef2}
Let $G$ be a locally compact group, and let $p \in [1,\infty)$. We say that $G$ has \emph{Reiter's property $(P_p)$} if there is a net $( m_\alpha )_\alpha$ of non-negative norm one functions in $L^p(G)$ such that
\[
  \lim_\alpha \| f L_\bullet(m_\alpha) - f \tensor m_\alpha \|_{{\cal C}_0(G) \tensor^\lambda L^p(G)} = 0
\]
for all $f \in {\cal C}_0(G)$.
\end{definition}
\section{Locally compact quantum groups---an overview} \label{LCQG}
In this section, we give a brief overview of locally compact quantum groups---as introduced by J.\ Kustermans and S.\ Vaes in \cite{KV1} and \cite{KV2}---with an emphasis on the von Neumann algebraic approach. For details, we refer to \cite{KV1}, \cite{KV2}, and \cite{vDa}. 
\par 
As a (von Neumann algebraic) locally compact quantum group is a Hopf--von Neumann algebra with additional structure, we begin with recalling the definition of a Hopf--von Neumann algebra ($\bar{\tensor}$ denotes the $W^\ast$-tensor product):
\begin{definition}
A \emph{Hopf--von Neumann algebra} is a pair $(\M,\Gamma)$, where $\M$ is a von Neumann algebra and $\Gamma \!: \M \to \M \bar{\tensor} \M$ is a \emph{co-multiplication}, i.e., a normal, unital, and injective $^\ast$-homomorphism satisfying $(\id \tensor \Gamma) \circ \Gamma = (\Gamma \tensor \id) \circ \Gamma$.
\end{definition}
\par
\begin{example}
For a locally compact group $G$, define $\Gamma_G \!: L^\infty(G) \to L^\infty(G \times G)$ by letting
\[
  (\Gamma_G \phi)(x,y) := \phi(xy) \qquad (\phi \in L^\infty(G), \, x,y \in G).
\]
Then $(L^\infty(G),\Gamma_G)$ is a Hopf--von Neumann algebra. 
\end{example}
\begin{remark}
Given a Hopf--von Neumann algebra $(\M,\Gamma)$, one can define a product $\ast$ on $\M_\ast$, the unique predual of $\M$, turning it into a Banach algebra:
\begin{equation} \label{prod}
  \langle f \ast g, x\rangle := \langle f \tensor g, \Gamma x\rangle \qquad (f,g \in \M_\ast, \,   x \in \M).
\end{equation}
If $G$ is a locally compact group, then applying (\ref{prod}) to $(L^\infty(G), \Gamma_G)$ yields the usual convolution product on $L^1(G)$.
\end{remark}
\par 
To define the additional structure that turns a Hopf--von Neumann algebra into a locally compact quantum group, we recall some basic facts about weights (see \cite{Tak2}, for instance). 
\par
Let $\M$ be a von Neumann algebra, and let $\M^+$ denote its positive elements. A \emph{weight} on $\M$ is an additive map $\phi \!: \M^+ \to [0,\infty]$ such that $\phi(tx) = t \phi(x)$ for $t \in [0,\infty)$ and $x \in \M^+$. We let
\[
  {\cal M}_\phi^+ := \{ x \in \M^+ : \phi(x) < \infty \}, \qquad
  {\cal M}_\phi := \text{the linear span of ${\cal M}_\phi^+$},
\]
and
\[
  {\cal N}_\phi := \{ x \in \M : x^\ast x \in {\cal M}_\phi \}.
\]
Then $\phi$ extends to a linear map on ${\cal M}_\phi$, and ${\cal N}_\phi$ is a left ideal of $\M$. Using the GNS-construction (\cite[p.\ 42]{Tak2}), we obtain a representation $\pi_\phi$ of $\M$ on some Hilbert space $\Hilbert_\phi$; we denote the canonical map from ${\cal N}_\phi$ into $\Hilbert_\phi$ by $\Lambda_\phi$. Moreover, we call $\phi$ \emph{ semifinite} if ${\cal M}_\phi$ is $w^\ast$-dense in $\M$, \emph{faithful} if $\phi(x) = 0$ for $x \in \M^+$ implies that $x = 0$, and \emph{normal} if $\sup_\alpha \phi(x_\alpha) = \phi\left( \sup_\alpha x_\alpha \right)$ for each bounded, increasing net $( x_\alpha )_\alpha$ in $\M^+$. If $\phi$ is faithful and normal, then the corresponding representation $\pi_\phi$ is faithful and normal, too (\cite[Proposition VII.1.4]{Tak2}).
\begin{definition} \label{lcqg}
A (von Neumann algebraic) \emph{locally compact quantum group} is a Hopf--von Neumann algebra $(\M,\Gamma)$ such that:
\begin{alphitems}
\item there is a normal, semifinite, faithful weight $\phi$ on $\M$---a \emph{left Haar weight}---which is left invariant, i.e., satisfies
\[
  \phi((f \tensor \id)(\Gamma x)) = \langle f, 1 \rangle \phi(x) \qquad   (f \in \M_\ast, \, x \in {\cal M}_\phi);
\]
\item there is a normal, semifinite, faithful weight $\psi$ on $\M$---a \emph{right Haar weight}---which is right invariant, i.e., satisfies
\[
  \psi((\id \tensor f)(\Gamma x)) = \langle f, 1 \rangle \psi(x) \qquad (f \in \M_\ast, \, x \in {\cal M}_\psi).
\]
\end{alphitems}
\end{definition}
\begin{example}
Let $G$ be a locally compact group. Then the Hopf--von Neumann algebra $(L^\infty(G),\Gamma_G)$ is a locally compact quantum group: $\phi$ and $\psi$ can be chosen as left and right Haar measure, respectively.
\end{example}
\begin{remarks}
\item Even though only the existence of a left and a right Haar weight, respectively, is presumed, both weights are actually unique up to a positive scalar multiple (see \cite{KV1} and \cite{KV2}). In order to make notation not too cumbersome, we shall thus simply write $(\M,\Gamma)$ for a locally compact quantum group whose left and right Haar weight will always be denoted by $\phi$ and $\psi$, respectively.
\item As discussed in \cite{KV1} and \cite{KV2}, locally compact quantum groups can equivalently be described in $\cstar$-algebraic terms. The $\cstar$-algebraic definition (\cite[Definition 4.1]{KV1}), however, is technically more involved, so that we shall not go into the details. 
\end{remarks}
\begin{definition} \label{Wdef}
Let $(\M,\Gamma)$ be a locally compact quantum group. The \emph{multiplicative unitary} of $(\M,\Gamma)$ is the unique operator $W \in
{\cal B}(\Hilbert_\phi \ttensor_2 \Hilbert_\phi)$, where $\ttensor_2$ stands for the Hilbert space tensor product, satisfying
\[
  W^\ast(\Lambda_\phi(x) \tensor \Lambda_\phi(y))
  = (\Lambda_\phi \tensor \Lambda_\phi)((\Gamma y)(x \tensor 1))
  \qquad (x,y \in {\cal N}_\phi).
\]
\end{definition}
\begin{example}
For a locally compact group $G$, the multiplicative unitary $W_G$ of $(L^\infty(G),\Gamma_G)$ is given by
\[
  (W_G\xi)(x,y) = \xi(x,x^{-1}y) \qquad (\xi \in L^2(G \times G), \, x ,y \in G).
\]
\end{example}
\begin{remarks}
\item Using the left invariance of $\phi$, it is easy to see that $W^\ast$ is an isometry whereas it is considerably more difficult to show that $W$ is indeed a unitary operator (\cite[Theorem 3.16]{KV1}).
\item The unitary $W$ lies in $\M \bar{\tensor} {\cal B}(\Hilbert_\phi)$ and implements the co-multiplication via 
\[
  \Gamma x = W^\ast (1 \tensor x) W \qquad (x \in \M)
\]
(see the discussion following \cite[Theorem 1.2]{KV2}). 
\item The definition of $W$ is made via the GNS-construction arising from $\phi$, so that one may want---in order to avoid confusion---rather speak of a \emph{left} multiplicative unitary. Indeed, one can define a \emph{right} multiplicative unitary in a similar fashion in terms of $\psi$: in \cite{JNR}, for instance, the right multiplicative unitary is used instead of the left one. It seems to be more or less a matter of taste with which of two multiplicative unitaries one prefers to work.
\end{remarks}
\par
To emphasize the parallels between locally compact quantum groups and groups, we shall use the following notation (which was suggested by Z.-J.\ Ruan and is also used in \cite{Run} and \cite{JNR}). We use the symbol $\G$ for a von Neumann algebraic, locally compact quantum group $(\M,\Gamma)$ and write: $L^\infty(\G)$ for $\M$, $L^1(\G)$ for $\M_\ast$, and $L^2(\G)$ for $\Hilbert_\phi$. If $L^\infty(\G) = L^\infty(G)$ for a locally compact group $G$ and $\Gamma = \Gamma_G$, we say that $\G$ actually \emph{is} a locally compact group, which is the case precisely if $L^\infty(\G)$ is abelian (this follows from \cite[Th\'eor\`eme 2.2]{BS}).
\par 
Given a locally compact quantum group $\G$ with multiplicative unitary $W$, we set
\[
  {\cal C}_0(\G) := \varcl{\{ (\id \tensor \omega)(W) : \omega \in {\cal B}(L^2(\G))_\ast \}}^{\| \cdot \|}
\]
It is relatively easy to see that ${\cal C}_0(\G)$ is a closed subalgebra of ${\cal B}(L^2(\G))$, but---which is much harder to show---it is even a $\cstar$-subalgebra. Restricting $\Gamma$ to ${\cal C}_0(\G)$ then yields a reduced $\cstar$-algebraic quantum group in the sense of \cite[Definition 4.1]{KV1} (see \cite[Proposition 1.6]{KV2}). If $\G$ is a locally compact group $G$, then ${\cal C}_0(G)$ just has the usual meaning: the continuous function on $G$ vanishing at infinity. Consequently, we write $M(\G)$ for ${\cal C}_0(\G)^\ast$. Like $L^1(\G)$, the dual space $M(\G)$ has a canonical product induced by $\Gamma$, turning it into a Banach algebra (\cite[p.\ 913]{KV1}) containing $L^1(\G)$ as a closed ideal (\cite[p.\ 914]{KV1}).
\par
Given a locally compact quantum group $\G$ with multiplicative unitary $W$, the \emph{left regular representation} of $\G$ is the map
\begin{equation} \label{leftreg}
  \lambda_2 \!: L^1(\G) \to {\cal B}(L^2(\G)), \quad f \mapsto (f \tensor \id)(W).
\end{equation}
Since $W \in L^\infty(\G) \bar{\tensor} {\cal B}(L^2(\G))$, it is clear that $\lambda_2$ is well defined, and it is easy to see that $\lambda_2$ is a contractive algebra homomorphism.
\begin{example}
For a locally compact group $G$, we have
\[
  (\lambda_2(f)\xi)(y) = \int_G f(x) \xi(x^{-1}y) \, dx \qquad (f \in L^1(G), \, \xi \in L^2(G)) 
\]
for almost all $y \in G$, i.e., $\lambda_2$ according to (\ref{leftreg}) is just the usual left regular representation of $L^1(G)$ on $L^2(G)$.
\end{example}
\par
Locally compact quantum groups allow for the development of a duality theory that extends Pontryagin duality for locally compact abelian groups.
\par
For a locally compact quantum group $\G$, set
\[
  L^\infty(\hat{\G}) : = 
  \varcl{\lambda_2(L^2(\G))}^{\,\text{$\sigma$-strongly$^\ast$}};
\]
it can be shown that $L^\infty(\hat{\G})$ is a von Neumann algebra. Let $\sigma$ denote the flip map on $L^2(\G) \ttensor_2 L^2(\G)$, i.e., $\sigma(\xi \tensor \eta) = \eta \tensor \xi$ for $\xi, \eta \in L^2(\G)$. Set $\hat{W} := \sigma W^\ast \sigma$. Then 
\[
  \hat{\Gamma} \!:  L^\infty(\hat{\G}) \to 
  L^\infty(\hat{\G}) \bar{\tensor}  L^\infty(\hat{\G}), \quad
  x \mapsto \hat{W}^\ast(1 \tensor x) \hat{W}
\]
is a co-multiplication. One can also define a left Haar weight $\hat{\phi}$ and a right Haar weight $\hat{\psi}$ for $(L^\infty(\hat{\G}),\hat{\Gamma})$ turning it into a locally compact quantum group again, the \emph{dual quantum group} of $\G$, which we denote by $\hat{\G}$, and whose multiplicative unitary is $\hat{W}$ as defined above. Finally, a Pontryagin duality theorem holds, i.e., $\Hat{\Hat{\G}} = \G$. For the details of this duality, we refer again to \cite{KV1} and \cite{KV2}.
\begin{example}
Let $G$ be a locally compact group. Since $L^\infty(\hat{G})$ is the $\sigma$-strong$^\ast$ closure of $\lambda_2(L^1(G))$, it equals $\VN(G)$, the group von Neumann algebra of $G$. Further, the co-multiplication $\hat{\Gamma}_G \!: \VN(G) \to \VN(G) \bar{\tensor} \VN(G)$ is given by
\[
  \hat{\Gamma}_G(\lambda(x)) = \lambda(x) \tensor \lambda(x) \qquad (x \in G). 
\]
Consequently, the product $\ast$ according to (\ref{prod}) on $\VN(G)_\ast$ is the usual pointwise product on $A(G)$, so that $L^1(\hat{G}) = A(G)$. The Plancherel weight on $\VN(G)$ (\cite[Definition VII.3.2]{Tak2}) is both a left and a right Haar weight for $(\VN(G),\hat{\Gamma}_G)$. Finally note that ${\cal C}_0(\hat{G})$ is the reduced group $\cstar$-algebra of $G$, so that $M(\hat{G})$ is the reduced Fourier--Stieltjes algebra from \cite{Eym}.
\end{example}
\section{$(P_1)$ for locally compact quantum groups} \label{P1}
With an eye on Definition \ref{Ppdef2}, we shall, in this section, formulate a version of property $(P_1)$ for locally compact quantum groups. To this end, we require the framework of operator space theory, as laid out in the monographs \cite{ER}, \cite{Pau}, and \cite{Pis}. We shall mostly follow \cite{ER} in our choice of notation; in particular, for two operator spaces $E$ and $F$, we denote the completely bounded operators from $E$ to $F$ by $\CB(E,F)$, we write $\| \cdot \|_\cb$ for the $\cb$-norm, and we use $\wTensor$ for the injective tensor product \emph{of operator spaces}. (Note that, if $\A$ and $\B$ are $\cstar$-algebras, then $\A \wTensor \B$ is just the spatial tensor product of $\cstar$-algebras.)
\par 
We begin with an elementary lemma:
\begin{lemma} \label{ellem}
Let $\Hilbert$ and $\Kilbert$ be Hilbert spaces, and let $A,B \in {\cal B}(\Kilbert) \wTensor {\cal K}(\Hilbert)$, where ${\cal K}(\Hilbert)$ denotes the compact operators on $\Hilbert$. Then the map
\begin{equation} \label{approxcb}
  {\cal B}(\Hilbert) \to {\cal B}(\Kilbert) \wTensor {\cal K}(\Hilbert), \quad
  x \mapsto A(1 \tensor x)B
\end{equation}
is completely bounded and belongs to the $\cb$-norm closure of the finite rank operators in $\CB({\cal B}(\Hilbert), {\cal B}(\Kilbert) \wTensor {\cal K}(\Hilbert))$.
\end{lemma}
\begin{proof}
The complete boundedness of (\ref{approxcb}) is clear.
\par 
To see that (\ref{approxcb}) is a norm limit of finite rank operators in $\CB({\cal B}(\Hilbert), {\cal B}(\Kilbert) \wTensor {\cal K}(\Hilbert))$, first note that it is enough to suppose that $A = S \tensor K$ and $B = T \tensor L$ with $S,T \in {\cal B}(\Kilbert)$ and $K, L \in {\cal K}(\Hilbert)$. Let $(K_n)_{n=1}^\infty$ and $( L_n )_{n=1}^\infty$ be finite rank operators on $\Hilbert$ such that
$K = \lim_{n \to \infty} K_n$ and $L = \lim_{n \to \infty} L_n$ in the norm topology of ${\cal B}(\Hilbert)$. For each $n \in \posints$, the operator
\[
  {\cal B}(\Hilbert) \to {\cal B}(\Kilbert) \wTensor {\cal K}(\Hilbert), \quad
  x \mapsto (S \tensor K_n)(1 \tensor x)(T \tensor L_n)
\]
has finite rank, and it is immediate that these operators converge to (\ref{approxcb}) in $\| \cdot \|_\cb$.
\end{proof}
\par 
Let $\G$ be a locally compact quantum group, and let $g \in L^1(\G)$. We define
\[
  (\Gamma|g) : L^\infty(\G) \to L^\infty(\G), \quad x \mapsto (\id \tensor g)(\Gamma x).
\]
It is immediate that $(\Gamma|g)$ is a weak$^\ast$-weak$^\ast$ continuous, completely bounded map.
\par 
For our next result---which will enable us to formulate property $(P_1)$ for locally compact quantum groups---, we use the following conventions:
\begin{itemize}
\item if $\A$ is any algebra, and $a$ and $b$ are any elements of $\A$, then $M_{a,b}$ denotes the two-sided multiplication map on $\A$ given by $M_{a,b}x := axb$ for $x \in \A$;
\item for any $\cstar$-algebra $\A$, its multiplier algebra (\cite[Definition III.6.22]{Tak1}) is denoted by ${\cal M}(\A)$;
\item if $\M$ is a von Neumann algebra on a Hilbert space $\Hilbert$ and $\xi$ and $\eta$ are vectors in $\Hilbert$, we write $\omega_{\xi,\eta}$ for the vector functional given by $\langle \omega_{\xi,\eta}, x \rangle = \langle x \xi,\eta \rangle$ for $x \in \M$.
\end{itemize}
We also recall that, if $E$ and $F$ are operator spaces, then the closure of the finite rank operators in $\CB(E,F)$ can be canonically identified with $F \wTensor E^\ast$ (\cite[Proposition 8.1.2]{ER}).
\begin{proposition} \label{approxprop}
Let $\G$ be a locally compact quantum group, let $g \in L^1(\G)$, and let $a,b \in {\cal C}_0(\G)$. Then $M_{a,b} \circ (\Gamma|g)$ is a completely bounded operator from $L^\infty(\G)$ to ${\cal C}_0(\G)$ that lies in the $\cb$-norm closure of the finite rank operators in $\CB(L^\infty(\G),{\cal C}_0(\G))$ and can be identified with an element of ${\cal C}_0(\G) \wTensor L^1(\G)$.
\end{proposition}
\begin{proof}
Since $L^\infty(\G)$ on $L^2(\G)$ is in standard form (\cite[Definition IX.1.13]{Tak2}), there are $\xi,\eta \in L^2(\G)$ such that $g = \omega_{\xi,\eta}$ (this follows from \cite[Lemma IX.1.6]{Tak2}). Choose $K,L \in {\cal K}(L^2(\G))$ such that $L \xi = \xi$ and $K^\ast \eta = \eta$ (clearly, rank one operators will do).
\par 
Let $W \in {\cal B}(L^2(\G) \ttensor_2 L^2(\G))$ be the multiplicative unitary of $\G$. By \cite[Proposition 3.21 and pp.\ 913--914]{KV1}---with the appropriate identifications made---, we have $W \in {\cal M}({\cal C}_0(\G) \wTensor {\cal K}(L^2(\G)))$, so that $(a \tensor K)W^\ast, W(b \tensor L) \in {\cal C}_0(\G) \wTensor {\cal K}(L^2(\G))$. From Lemma \ref{ellem}, we conclude that
\begin{equation} \label{approxcb2}
  L^\infty(\G) \to {\cal B}(L^2(\G)) \wTensor {\cal K}(L^2(\G)), \quad x \mapsto (a \tensor K)W^\ast (1 \tensor x) W(b \tensor L)
\end{equation}
is a norm limit of finite rank operators in $\CB(L^\infty(\G), {\cal B}(L^2(\G)) \wTensor {\cal K}(L^2(\G)))$. It is straightforward that (\ref{approxcb2}) actually lies in $\CB(L^\infty(\G), {\cal C}_0(\G) \wTensor {\cal K}(L^2(\G)))$.
\par 
Let $\xi', \eta' \in L^2(\G)$, and note that, for $x \in L^\infty(\G)$, we have
\[
  \begin{split}
  \langle ((M_{a,b} \circ (\Gamma|g))x)\xi', \eta' \rangle 
  & = \langle a (\id \tensor g)(W^\ast (1 \tensor x) W) b \xi', \eta' \rangle \\
  & = \langle (a \tensor 1)W^\ast(1 \tensor x) W(b \tensor 1) (\xi' \tensor \xi), \eta' \tensor \eta \rangle \\
  & = \langle (a \tensor 1)W^\ast(1 \tensor x) W(b \tensor 1) (\xi' \tensor L \xi), \eta' \tensor K^\ast \eta \rangle \\
  & = \langle (a \tensor K)W^\ast(1 \tensor x) W(b \tensor L) (\xi' \tensor \xi), \eta' \tensor \eta \rangle \\
  & = \langle (\id \tensor g)((a \tensor K)W^\ast(1 \tensor x) W(b \tensor L)) \xi', \eta' \rangle.
  \end{split}
\]
Since $\xi', \eta' \in L^2(\G)$ were arbitrary, this means that
\[
  (M_{a,b} \circ (\Gamma|g))x = (\id \tensor g)((a \tensor K)W^\ast(1 \tensor x) W(b \tensor L))  \qquad (x \in L^\infty(\G)),
\]
i.e., $M_{a,b} \circ (\Gamma|g)$ is the composition of (\ref{approxcb2}) with the Tomiyama slice map $\id \tensor g$ and thus is a norm limit of finite rank operators in $\CB(L^\infty(\G),{\cal C}_0(\G))$.
\par 
By \cite[Proposition 8.1.2]{ER}, $M_{a,b} \circ (\Gamma|g)$ can be canonically identified with an element of ${\cal C}_0(\G) \wTensor L^\infty(\G)^\ast$. In order to prove that $M_{a,b} \circ (\Gamma|g)$ actually lies in ${\cal C}_0(\G) \wTensor L^1(\G)$, we show that $(M_{a,b} \circ (\Gamma|g))^\ast \!: M(\G) \to L^\infty(\G)^\ast$ attains its values in $L^1(\G)$. For $\mu \in M(\G)$ and $x \in L^\infty(\G)$, we have
\[
  \begin{split}
  \langle (M_{a,b} \circ (\Gamma|g))^\ast \mu, x \rangle & = \langle (M_{a,b} \circ (\Gamma|g)) x, \mu \rangle \\
  & = \langle (a \tensor 1)(\Gamma x)(b \tensor 1), \mu \tensor g \rangle \\
  & = \langle \Gamma x, b\mu a \tensor g \rangle \\
  & = \langle b \mu a \ast g, x \rangle, 
  \end{split}
\]
so that
\begin{equation} \label{adjoint}
  (M_{a,b} \circ (\Gamma|g))^\ast \mu = b\mu a \ast g.
\end{equation}
(We denote the canonical module actions of a $\cstar$-algebra on its dual by juxtaposition.) Since $L^1(\G)$ is an ideal in $M(\G)$, it follows from (\ref{adjoint}) that $(M_{a,b} \circ (\Gamma|g))^\ast M(\G) \subset L^1(\G)$, so that $M_{a,b} \circ (\Gamma|g)$ is canonically represented by an element of ${\cal C}_0(\G) \wTensor L^1(\G)$.
\end{proof}
\par 
Let $G$ be a locally compact group, let $g \in L^1(G)$, and let $a,b \in {\cal C}_0(G)$. Then $M_{a,b} \circ (\Gamma|g) \in {\cal C}_0(G) \wTensor L^1(G) = {\cal C}_0(G) \tensor^\lambda L^1(G)$ is nothing but $ab L_\bullet(g)$ in the notation of Section \ref{prelude}, as a routine verification shows.
\par 
With Definition \ref{Ppdef2} in mind, we can thus extend property $(P_1)$ from locally compact groups to locally compact quantum groups:
\begin{definition} \label{P1def}
A locally compact quantum group $\G$ is said to have \emph{Reiter's property $(P_1)$} if there is a net $( m_\alpha )_\alpha$ of states in $L^1(\G)$ such that
\[
  \lim_\alpha \| M_{a,b} \circ (\Gamma|m_\alpha) - ab \tensor m_\alpha \|_{{\cal C}_0(\G) \wTensor L^1(\G)} = 0
\] 
for all $a,b \in {\cal C}_0(\G)$.
\end{definition}
\section{Amenability and $(P_1)$}
Recall the definition of an amenable, locally compact quantum group:
\begin{definition} \label{amdef}
A locally compact quantum group $\G$ is called \emph{amenable} if it has a \emph{left invariant mean}, i.e., a state $M$ on $L^\infty(\G)$ such that
\begin{equation} \label{LIM}
  \langle (f \tensor \id)(\Gamma x), M \rangle = \langle f, 1 \rangle \langle x, M \rangle \qquad (f \in L^1(\G), \,
  x \in L^\infty(\G)).
\end{equation}
\end{definition}
\begin{remarks}
\item Our use of the term amenable is the same as in \cite{BT}, but there is no general consensus in the literature about terminology: an amenable, locally compact quantum group according to Definition \ref{amdef} is called \emph{Voiculescu amenable} in \cite{Rua2} and \emph{weakly amenable} in \cite{DQV}.
\item There is an element of asymmetry in Definition \ref{amdef}: a state $M$ on $L^\infty(\G)$ is called a \emph{right invariant mean} if
\begin{equation} \label{RIM}
  \langle (\id  \tensor f)(\Gamma x), M \rangle = \langle f, 1 \rangle \langle x, M \rangle \qquad (f \in L^1(\G), \,
  x \in L^\infty(\G))
\end{equation}
holds and an \emph{invariant mean} if both (\ref{LIM}) and (\ref{RIM}) are satisfied. So, $\G$ is amenable if and only if there is a \emph{left} invariant mean on $L^\infty(\G)$. However, by \cite[Proposition 3]{DQV}, the amenability of $\G$ already implies the existence of an invariant mean.
\item The standard approximation argument (see \cite{ES1}, for instance) immediately yields that $\G$ is amenable if and only if there is a net $( m_\alpha )_\alpha$ of states in $L^1(\G)$ such that
\begin{equation} \label{asyLIM}
  \lim_\alpha \| f \ast m_\alpha - \langle f, 1 \rangle m_\alpha \| = 0 \qquad (f \in L^1(\G)).
\end{equation}
\item If $G$ is a locally compact group, then a state $M$ as in Definition \ref{amdef} is topologically left invariant in the sense of \cite[Definition 4.3]{Pie}. By \cite[Theorem 4.19]{Pie}, this means that $G$ is amenable in the sense of Definition \ref{amdef} if and only if it is amenable in the classical sense.
\end{remarks}
\par
It is easy to see that $(P_1)$ implies amenability:
\begin{proposition} \label{P1toam}
Let $\G$ be a locally compact quantum group with Reiter's property $(P_1)$. Then $\G$ is amenable.
\end{proposition}
\begin{proof}
Let $( m_\alpha )_\alpha$ be a net as in Definition \ref{P1def}, and let $f \in L^1(\G)$. By Cohen's factorization theorem (\cite[Corollary 2.9.26]{Dal}), there are $a,b \in {\cal C}_0(\G)$ and $g \in L^1(\G)$ such that $f = bga$. 
\par 
For any Banach space $E$, we denote its closed unit ball by $\Ball(E)$. 
\par 
We then have:
\[
  \begin{split}
  \| f \ast m_\alpha - \langle 1,f \rangle m_\alpha \|  
  & = \sup \{ | \langle f \tensor m_\alpha, \Gamma x - 1 \tensor x \rangle | : x \in \Ball(L^\infty(\G)) \} \\
  & = \sup \{ | \langle bga \tensor m_\alpha, \Gamma x - 1 \tensor x \rangle | : x \in \Ball(L^\infty(\G)) \} \\
  & = \sup \{ | \langle g \tensor m_\alpha, (a \tensor 1)(\Gamma x)(b\tensor 1) - ab \tensor x \rangle | 
    : x \in \Ball(L^\infty(\G)) \} \\
  & = \sup \{ | \langle g, (M_{a,b} \circ (\Gamma|m_\alpha))x - \langle m_\alpha, x \rangle ag \rangle | : 
                x \in \Ball(L^\infty(\G)) \} \\
  & \leq \| g \| \sup \{ \| (M_{a,b} \circ (\Gamma|m_\alpha))x - \langle m_\alpha, x \rangle ab\| : x \in \Ball(L^\infty(\G)) \} \\
  & \leq \| g \| \| M_{a,b} \circ (\Gamma|m_\alpha) - ab \tensor m_\alpha \|_{{\cal C}_0(\G) \wTensor L^1(\G)} \\
  & \stackrel{\alpha} {\longrightarrow} 0.
  \end{split}
\]
It is clear that any weak$^\ast$ accumulation point of $(m_\alpha)_\alpha$ in $L^\infty(\G)^\ast$ is a left invariant mean.
\end{proof}
\par 
For the converse of Proposition \ref{P1toam}, we require a few preparations.
\par 
Let $E$ be an operator space; deviating from \cite{ER}, but for the sake of notational clarity, we denote, for $n \in \posints$, the $n$-th matrix level of $E$ by $\mathbb{M}_n(E)$. A \emph{matricial subset} of $E$ is a sequence $\boldsymbol{S} = (S_n )_{n=1}^\infty$ with $S_n \subset \mathbb{M}_n(E)$ for $n \in \posints$. We use the usual set theoretic symbols for matricial points and subsets termwise, e.g., if $\boldsymbol{S} = (S_n)_{n=1}^\infty$ and $\boldsymbol{T} = ( T_n )_{n=1}^\infty$ are matricial subsets of $E$, then $\boldsymbol{S} \cup \boldsymbol{T}$ is defined as $( S_n \cup T_n )_{n=1}^\infty$.
\par 
Given two operator spaces $E$ and $F$, $n \in \posints$, and a linear map $T \!: E \to F$, we write (again, not following \cite{ER}) $T^{(n)} \!: \mathbb{M}_n(E) \to \mathbb{M}_n(F)$ for the $n$-th amplification of $T$.
\begin{definition} \label{cuniform}
Let $E$ and $F$ be operator spaces, let $( T_\alpha )_\alpha$ be a net in $\CB(E,F)$, let $T \in \CB(E,F)$, and let $\boldsymbol{S} = (S_n)_{n=1}^\infty$ be a matricial subset of $E$. We say that $(T_\alpha)_\alpha$ converges to $T$ \emph{completely uniformly on $\boldsymbol{S}$} if
\[
  \lim_\alpha \sup_{n \in \posints} \sup \left\{ \left\| T^{(n)}_\alpha  x - T^{(n)} x \right\|_n : x \in S_n \right\} \to 0.
\]
\end{definition}
\begin{lemma} \label{cunilem}
Let $E_1, \ldots, E_m$, $E$, and $F$ be operator spaces, and let $S_j \in \CB(E_j,E)$ for $j =1, \ldots, m$ lie in the $\cb$-norm closure of the finite rank operators. Let 
\[
  \boldsymbol{K}_j := \left( S_j^{(n)} (\Ball(\mathbb{M}_n(E_j))) \right)_{n=1}^\infty \qquad (j =1, \ldots, m),
\]
and set $\boldsymbol{K} := \boldsymbol{K}_1 \cup \cdots \cup \boldsymbol{K}_m$. Then every norm bounded net $( T_\alpha )_\alpha$ in $\CB(E,F)$ that converges to $T \in \CB(E,F)$ pointwise on $E$ converges to $T$ completely uniformly on $\boldsymbol{K}$.
\end{lemma}
\begin{proof}
Suppose without loss of generality that $m = 1$.
\par 
The completely uniform convergence of $( T_\alpha )_\alpha$ to $T$ on $\boldsymbol{K}_1$ amounts to $\| T_\alpha S_1 - T S_1 \|_\cb \to 0$. Since $( T_\alpha )_\alpha$ is norm bounded in $\CB(E,F)$ and since $S_1$ is a norm limit of finite rank operators in $\CB(E_1,E)$, there is no loss of generality to suppose that $S_1$ \emph{is} a finite rank operator. 
\par 
Let $E_0$ be a finite-dimensional subspace of $E$ with $S_1 E_1 \subset E_0$. Since $\dim E_0 < \infty$, the identity map $\id_{E_0} \!: E_0 \to \max E_0$ is completely bounded (\cite[Theorem 14.3(ii)]{Pau}). Hence, we have the (Banach space) isomorphisms
\[
  \CB(E_0,F) \cong \CB(\max E_0,F) \cong {\cal B}(E_0,F),
\]
where the last isomorphism holds by the definition of $\max E$ and is, in fact, isometric (\cite[(3.3.9)]{ER}). Since the unit ball of $E_0$ is compact, and since $( T_\alpha )_\alpha$ is norm bounded in ${\cal B}(E_0,F)$, we conclude that $T_\alpha |_{E_0} \to T |_{E_0}$ in the norm topology of ${\cal B}(E_0,F)$ and thus of $\CB(E_0,F)$. Finally, note that
\[
  \| T_\alpha S_1 - T S_1 \|_\cb \leq \frac{1}{\max\{ \|S_1 \|_\cb, 1 \} } \| T_\alpha |_{E_0} - T |_{E_0} \|_\cb \to 0,
\]
which completes the proof.
\end{proof}
\begin{remark}
Let $E$ and $F$ be Banach spaces, let $( T_\alpha )_\alpha$ be a norm bounded net in ${\cal B}(E,F)$, and let $T \in {\cal B}(E,F)$ be such that $T_\alpha \to T$ pointwise on $E$. Then it is elementary (and was used in the proof of Lemma \ref{cunilem}) that $T_\alpha \to T$ uniformly on all compact subsets of $E$. Lemma \ref{cunilem} is a fairly crude attempt to adapt this fact to an operator space setting. One major obstacle to establishing a more satsifactory operator space variant is the apparent difficulty of finding a proper notion of compactness suited for operator spaces (see \cite{Web} and \cite{Yew}).
\end{remark}
\par 
We can now prove the first main result of this paper:
\begin{theorem} \label{P1thm}
Let $\G$ be a locally compact quantum group. Then the following are equivalent:
\begin{items}
\item $\G$ is amenable;
\item $\G$ has Reiter's property $(P_1)$.
\end{items}
\end{theorem}
\begin{proof}
As (ii) $\Longrightarrow$ (i) is Proposition \ref{P1toam}, all we need to prove is (i) $\Longrightarrow$ (ii).
\par 
Let $a_1, b_1, \ldots, a_\nu, b_\nu \in {\cal C}_0(\G)$, and let $\epsilon > 0$. We need to show that there is a state $m \in L^1(\G)$ such that
\begin{equation} \label{goal}
  \| M_{a_j,b_j} \circ (\Gamma|m) - a_j b_j \tensor m \|_{{\cal C}_0(\G) \wTensor L^1(\G)} < \epsilon \qquad (j=1,\ldots, \nu).
\end{equation}
\par 
Since $\G$ is amenable, there is a net $( m_\alpha )_{\alpha \in \mathbb A}$ of states in $L^1(\G)$ such that (\ref{asyLIM}) holds.
For $\alpha \in \mathbb A$, define
\[
  T_\alpha \!: L^1(\G) \to L^1(\G), \quad f \mapsto f \ast m_\alpha - \langle f, 1 \rangle m_\alpha.
\]
The net $( T_\alpha )_\alpha$ lies in $\CB(L^1(\G))$, is norm bounded, and converges to $0$ pointwise on $L^1(\G)$ by (\ref{asyLIM}). 
\par 
Let $m_0 \in L^1(\G)$ be an arbitrary state. For $j = 1, \ldots, \nu$, let the matricial subset $\boldsymbol{K}_j = ( K_{j,n} )_{n=1}^\infty$ of $L^1(\G)$ be defined through
\[
  K_{j,n} := \{ [ b_j \mu_{k,l} a_j \ast m_0 ] : [ \mu_{k,l} ] \in \Ball(\mathbb{M}_n(M(\G))) \} \qquad (n \in \posints).
\]
For $j = 1, \ldots, \nu$, let $S_j \in \CB(M(\G),L^1(\G))$ be defined as $(M_{a_j,b_j} \circ (\Gamma|m_0))^\ast$. By Proposition \ref{approxprop}, this means that $S_j$ belongs to the norm closure of the finite rank operators in $\CB(M(\G),L^1(\G))$. A simple calculation shows that
\[
  S_j \mu = b_j \mu a_j \ast m_0 \qquad ( j =1, \ldots,\nu, \, \mu \in M(\G)),
\]
so that
\[
  K_{j,n} = S_j^{(n)} ( \Ball(\mathbb{M}_n(M(\G))) ) \qquad (j = 1, \ldots, \nu, \, n \in \posints).
\]
Invoking Lemma \ref{cunilem}---with $\boldsymbol{K}_1, \ldots, \boldsymbol{K}_\nu$ as just defined---as well as (\ref{asyLIM}), we obtain $\alpha_\epsilon \in \mathbb A$ such that
\[
  \sup_{n \in \posints} 
  \sup \left\{ \left\| T_{\alpha_\epsilon}^{(n)} f \right\| : f \in K_{1,n} \cup \cdots \cup K_{\nu,n} \right\} < \frac{\epsilon}{2} 
\]
as well as
\[
  \| m_0 \ast m_{\alpha_\epsilon} - m_{\alpha_\epsilon} \| 
  < \frac{1}{\max \{ \| a_1 b_1 \|, \ldots, \| a_\nu b_\nu \|, 1 \}} \frac{\epsilon}{2}.
\]
Set $m := m_0 \ast m_{\alpha_\epsilon}$.
\par 
To see that (\ref{goal}) holds, first note that
\[
  \begin{split}
  \| M_{a_j,b_j} \circ (\Gamma|m) - a_j b_j \tensor m \| 
  & \leq \| M_{a_j,b_j} \circ (\Gamma|m) - a_j b_j \tensor m_{\alpha_\epsilon} \| 
    + \| a_j b_j \tensor m_{\alpha_\epsilon} - a_j b_j \tensor m \| \\
  & = \| M_{a_j,b_j} \circ (\Gamma|m) - a_j b_j \tensor m_{\alpha_\epsilon} \| 
    + \| a_j b_j \| \| m - m_{\alpha_\epsilon} \| \\
  & < \| M_{a_j,b_j} \circ (\Gamma|m) - a_j b_j \tensor m_{\alpha_\epsilon} \| + \frac{\epsilon}{2} 
  \qquad ( j =1, \ldots, \nu).
\end{split}
\]
In order to establish (\ref{goal}), it is thus enough to show that 
\begin{equation} \label{subgoal}
  \| M_{a_j,b_j} \circ (\Gamma|m) - a_j b_j \tensor m_{\alpha_\epsilon} \| < \frac{\epsilon}{2} \qquad (j=1,\ldots, \nu).
\end{equation}
\par 
With $j \in \{1, \ldots, \nu \}$ fixed, note that
\begin{multline} \label{sups}
  \| M_{a_j,b_j} \circ (\Gamma|m) - a_j b_j \tensor m_{\alpha_\epsilon} \| \\
  = \sup_{n \in \posints} \sup \{ \| [ (M_{a_j,b_j} \circ (\Gamma|m)) x_{k,l}  
      - a_jb_j \langle m_{\alpha_\epsilon}, x_{k,l} \rangle ] \|_n : [ x_{k,l} ] \in \Ball(\mathbb{M}_n(L^\infty(\G))) \}.
\end{multline}
Let $\llangle \cdot, \cdot \rrangle$ denote the matrix duality of \cite{ER}. By \cite[(3.2.4)]{ER}, the second supremum of the right hand side of (\ref{sups}) is then computed as
\begin{multline} \label{sup2}
  \sup \{ \| \llangle  [ (M_{a_j,b_j} \circ (\Gamma|m))x_{k,l} 
      - \langle m_{\alpha_\epsilon}, x_{k,l} \rangle a_j b_j  ], [\mu_{\kappa,\lambda}] \rrangle \|_{n^2}: \\
     [ x_{k,l} ] \in \Ball(\mathbb{M}_n(L^\infty(\G))), \, [ \mu_{\kappa,\lambda} ] \in \Ball(\mathbb{M}_n(M(\G))) \}.
\end{multline}
For $x \in L^\infty(\G)$ and $\mu \in M(\G)$, we have
\[
  \begin{split}
  \lefteqn{\langle (M_{a_j,b_j} \circ (\Gamma|m)) x- \langle m_{\alpha_\epsilon}, x \rangle a_j b_j, \mu \rangle} \\
  & = \langle \Gamma x , b_j \mu a_j \tensor m \rangle - \langle 1 \tensor x, b_j \mu a_j \tensor m_{\alpha_\epsilon} \rangle \\
  & = \langle \Gamma x , b_j \mu a_j \tensor (m_0 \ast m_{\alpha_\epsilon}) \rangle 
    - \langle 1 \tensor x, b_j \mu a_j \tensor m_{\alpha_\epsilon} \rangle \\
  & = \langle (\id \tensor \Gamma)(\Gamma x) , b_j \mu a_j \tensor m_0 \tensor m_{\alpha_\epsilon} \rangle 
    - \langle 1 \tensor x, b_j \mu a_j \tensor m_{\alpha_\epsilon} \rangle \\
  & = \langle (\Gamma \tensor \id)(\Gamma x) , b_j \mu a_j \tensor m_0 \tensor m_{\alpha_\epsilon} \rangle 
    - \langle 1 \tensor x, b_j \mu a_j \tensor m_{\alpha_\epsilon} \rangle \\
  & = \langle (\Gamma \tensor \id)(\Gamma x) - \Gamma 1 \tensor x, b_j \mu a_j \tensor m_0 \tensor m_{\alpha_\epsilon} \rangle \\
  & = \langle \Gamma x - 1 \tensor x, (b_j \mu a_j \ast m_0) \tensor m_{\alpha_\epsilon} \rangle \\
  & = \langle (b_j \mu a_j \ast m_0) \ast m_{\alpha_\epsilon} - \langle b_j \mu a_j \ast m_0, 1 \rangle m_{\alpha,\epsilon}, x \rangle.
  \end{split}
\]
Again from \cite[(3.2.4)]{ER}, we therefore conclude that (\ref{sup2}) equals
\[
  \sup \{ \| [(b_j \mu_{\kappa,\lambda} a_j \ast m_0) \ast m_{\alpha_\epsilon} 
  - \langle b_j \mu_{\kappa,\lambda} a_j \ast m_0, 1 \rangle m_{\alpha_\epsilon} ] \|_n : 
  [ \mu_{\kappa,\lambda} ] \in \Ball(\mathbb{M}_n(M(\G))) \}.
\]
We thus have 
\[
  \begin{split}
  \lefteqn{\| M_{a_j,b_j} \circ (\Gamma|m) - a_j b_j \tensor m_{\alpha_\epsilon} \|} & \\ 
  & = \sup_{n \in \posints} \\
  & \qquad \sup \{ \|  [(b_j \mu_{\kappa,\lambda} a_j \ast m_0) \ast m_{\alpha_\epsilon} 
    - \langle b_j \mu_{\kappa,\lambda} a_j \ast m_0, 1 \rangle m_{\alpha_\epsilon} ] \|_n : [\mu_{\kappa,\lambda} ] \in \Ball(\mathbb{M}_n(M(\G))) \} \\
  & \leq \sup_{n \in \posints} \sup \{  \| [ f_{\kappa,\lambda} \ast m_{\alpha_\epsilon} 
    - \langle f_{\kappa,\lambda}, 1 \rangle m_{\alpha_\epsilon} ] \|_n : 
    [ f_{\kappa,\lambda} ] \in K_{j,n} \} \\
  & = \sup_{n \in \posints} \sup \left\{ \| [ T_{\alpha_\epsilon} f_{\kappa,\lambda} ] \|_n : 
    [ f_{\kappa,\lambda} ] \in K_{j,n} \right\} \\
  & < \frac{\epsilon}{2}, \qquad\text{by the choice of $\alpha_\epsilon$},
  \end{split}
\]
for $j =1, \ldots, \nu$, i.e., (\ref{subgoal}) holds.
\end{proof}
\section{$(P_2)$ and co-amenability}
We finally turn to defining Reiter's property $(P_2)$ for locally compact quantum groups.
\par
Following \cite{ER}, we denote the column and row operator space over a Hilbert space $\Hilbert$ by $\Hilbert_c$ and $\Hilbert_r$, respectively. Given $T \in {\cal B}(\Hilbert) \bar{\tensor} {\cal B}(\Hilbert)$ and $\xi \in \Hilbert$, we have a linear map
\[
  (T|\xi) : \overline{\Hilbert} \to {\cal B}(\Hilbert), \quad \eta \mapsto (\id \tensor \omega_{\xi,\eta})(T),
\]
where $\overline{\Hilbert}$ denotes the complex conjugate Hilbert space of $\Hilbert$. From the definition of row Hilbert space, it is routine to verify that $(T|\xi) \in \CB\left(\overline{\Hilbert}_r,{\cal B}(\Hilbert) \right)$. By \cite[p.\ 59]{ER}, we can canonically identity $\overline{\Hilbert}_r$ with $\Hilbert_c^\ast$, and thus view $(T|\xi)$ as an operator in $\CB(\Hilbert_c^\ast, {\cal B}(\Hilbert))$.
\par 
We have the following $L^2$-analog of Proposition \ref{approxprop}:
\begin{proposition} \label{approxprop2}
Let $\G$ be a locally compact quantum group with multiplicative unitary $W$, let $\xi \in L^2(\G)$, and let $a,b \in {\cal C}_0(\G)$. Then $M_{a,b} \circ (W|\xi)$ is a completely bounded operator from $L^2(\G)_c^\ast$ to ${\cal C}_0(\G)$ that lies in the $\cb$-norm closure of the finite rank operators in $\CB(L^2(\G)^\ast_c,{\cal C}_0(\G))$ and can be identified with an element of ${\cal C}_0(\G) \wTensor L^2(\G)_c$.
\end{proposition}
\begin{proof}
Choose $L \in {\cal K}(L^2(\G))$ with $L \xi = \xi$, so that $(a \tensor 1)W(b \tensor L) \in {\cal C}_0(\G) \wTensor {\cal K}(L^2(\G))$. By the definition of $L^2(\G)_c$, the linear map
\[
  T_\xi \!: {\cal K}(L^2(\G)) \to L^2(\G)_c, \quad K \mapsto K\xi
\]
is completely bounded, so that
\[
  (\id \tensor T_\xi)((a \tensor 1)W(b \tensor L)) \in {\cal C}_0(\G) \wTensor L^2(\G)_c.
\]
Embedding ${\cal C}_0(\G) \wTensor L^2(\G)_c$ canonically into $\CB(L^2(\G)_c^\ast, {\cal C}_0(\G))$, we see that
\[
  \begin{split}
  (\id \tensor T_\xi)((a \tensor 1)W(b \tensor L))\eta & = (\id \tensor \omega_{\xi,\eta})((a \tensor 1)W(b \tensor L)) \\
  & = M_{a,b} \circ (\id \tensor \omega_{L\xi,\eta})(W) \\
  & = M_{a,b} \circ (\id \tensor \omega_{\xi,\eta})(W) \\
  & = (M_{a,b} \circ (W|\xi))\eta \qquad (\eta \in \Hilbert),
  \end{split}
\]
which completes the proof.
\end{proof}
\par
Let $G$ be a locally compact group, let $a,b \in {\cal C}_0(G)$, and let $\xi \in L^2(G)$. Then it is easily checked that $M_{a,b} \circ (W|\xi) = ab L_\bullet(\xi)$. With this in mind, we define:
\begin{definition} \label{P2def}
Let $\G$ be a locally compact quantum group with multiplicative unitary $W$. We say that $\G$ has \emph{Reiter's property $(P_2)$} if there is a net $( \xi_\alpha )_\alpha$ of unit vectors in $L^2(\G)$ such that
\[
  \lim_\alpha \| M_{a,b} \circ (W|\xi_\alpha) - ab \tensor \xi_\alpha \|_{{\cal C}_0(\G) \wTensor L^2(\G)_c} = 0
\] 
for all $a,b \in {\cal C}_0(\G)$.
\end{definition}
\begin{remarks}
\item Let $G$ be a locally compact group with $(P_2)$ in the sense of Definition \ref{Ppdef2}, and let $(\xi_\alpha )_\alpha$ be a net in $L^2(G)$ as required by that definition; then $( \xi_\alpha)_\alpha$ clearly satisfies Definition \ref{P2def}. On the other hand, if $G$ has property $(P_2)$ in the sense of Definition \ref{P2def} and if $(\xi_\alpha)_\alpha$ is a corresponding net of unit vectors in $L^2(G)$, then $(  | \xi_\alpha |)_\alpha$, where the modulus is taken pointwise almost everywhere, satisfies Definition \ref{Ppdef}.
\item Since $L^\infty(\G)$ is in standard form on $L^2(\G)$ (\cite[Definition IX.1.13]{Tak2}), there is a canonical self-dual cone $L^2(\G)_+$ in $L^2(\G)$ that provides a notion of positivity in $L^2(\G)$. We could thus have required the net $(\xi_\alpha)_\alpha$ in Definition \ref{P2def} to be from $L^2(\G)_+$. The reason why we haven't done this is Theorem \ref{P2thm} below: we do not know if it remains true with this additional requirement. (Unlike in the group case, we cannot conclude that, if a net as in Definition \ref{P2def} exists, then it can always be found in $L^2(\G)_+$; see also the remark immediately after the proof of Theorem \ref{P2thm}.)
\end{remarks}
\par 
For our second main result, recall the definition of a co-amenable, locally compact quantum group:
\begin{definition} \label{coamdef}
A locally compact quantum group $\G$ is called \emph{co-amenable} if the Banach algebra $L^1(\G)$ has a bounded approximate identity.
\end{definition}
\begin{remarks}
\item There are several descriptions of co-amenability equivalent to Definition \ref{coamdef}: see \cite[Theorem 3.1]{BT}. In particular, $\G$ is co-amenable if and only if there is a net $( \xi_\alpha )_\alpha$ of unit vectors in $L^2(\G)$ such that
\[
   \| W(\xi_\alpha \tensor \eta) - \xi_\alpha \tensor \eta \| \to 0 \qquad (\eta \in L^2(\G)).
\]
\item If $\hat{\G}$ is co-amenable (\cite[Theorem 3.2]{BT}), then $\G$ is amenable whereas the converse is unknown unless $\G$ is discrete (\cite{Tom}) or a group (\cite{Lep}).
\end{remarks}
\begin{theorem} \label{P2thm}
Let $\G$ be a locally compact quantum group. Then the following are equivalent:
\begin{items}
\item $\G$ has $(P_2)$;
\item $\hat{\G}$ is co-amenable.
\end{items}
\end{theorem}
\begin{proof}
(i) $\Longrightarrow$ (ii): Let $(\xi_\alpha )_\alpha$ be a net as required by Definition \ref{P2def}. We claim that
\[
  \| \hat{W}(\xi_\alpha \tensor \eta) - \xi_\alpha \tensor \eta \| \to 0 \qquad (\eta \in L^2(\G))
\]
or rather---equivalently by the definition of $\hat{W}$---
\begin{equation} \label{asycoam}
  \| W (\eta \tensor \xi_\alpha) - \eta \tensor \xi_\alpha \| \to 0 \qquad (\eta \in L^2(\G)).
\end{equation}
\par 
Let $\eta \in L^2(\G)$ be a unit vector, and use Cohen's factorization theorem (\cite[Corollary 2.9.26]{Dal}) to obtain $a \in {\cal C}_0(\G)$ and $\zeta \in L^2(\G)$ such that $\eta = a \zeta$. By Definition \ref{P2def}, 
\[
  \| M_{a^\ast, a} \circ (W|\xi_\alpha) - a^\ast a \tensor \xi_\alpha \|_{{\cal C}_0(\G) \wTensor L^2(\G)_c} \to 0
\]
holds. By the definition of column Hilbert space, the map
\[
  T_\zeta \!: {\cal C}_0(\G) \to L^2(\G)_c, \quad b \mapsto b\zeta
\]
lies in $\CB({\cal C}_0(\G), L^2(\G)_c)$ with $\| T_\zeta \|_\cb \leq \| \zeta \|$. Since $L^2(\G)_c \wTensor L^2(\G)_c = (L^2(\G) \ttensor_2 L^2(\G))_c$ (\cite[Proposition 9.3.5]{ER}), it follows that
\[
  \begin{split}
  \lefteqn{\| (a^\ast \tensor 1)W(\eta \tensor \xi_\alpha) - a^\ast \eta \tensor \xi_\alpha \|_{L^2(\G) \ttensor_2 L^2(\G)}} & \\
  & = \| (a^\ast \tensor 1)W(\eta \tensor \xi_\alpha) - a^\ast \eta \tensor \xi_\alpha \|_{L^2(\G)_c \wTensor L^2(\G)_c} \\
  & = \| (T_\zeta \tensor \id)(M_{a^\ast, a} \circ (W|\xi_\alpha) - a^\ast a \tensor \xi_\alpha) \|_{L^2(\G)_c \wTensor L^2(\G)_c} \\
  & \leq \| \zeta \| \| M_{a^\ast, a} \circ (W|\xi_\alpha) - a^\ast a \tensor \xi_\alpha \|_{{\cal C}_0(\G) \wTensor L^2(\G)_c} \\
  & \to 0
  \end{split}
\]
and thus
\[
  \begin{split}
  1 & = \lim_\alpha \langle \eta \tensor \xi_\alpha, \eta \tensor \xi_\alpha \rangle \\
  & = \lim_\alpha \langle a \zeta \tensor \xi_\alpha, \eta \tensor \xi_\alpha \rangle \\
  & = \lim_\alpha \langle \zeta \tensor \xi_\alpha, a^\ast\eta \tensor \xi_\alpha \rangle \\
  & = \lim_\alpha \langle \zeta \tensor \xi_\alpha, (a^\ast \tensor 1)W(\eta \tensor \xi_\alpha) \rangle \\
  & = \lim_\alpha \langle a\zeta \tensor \xi_\alpha, W(\eta \tensor \xi_\alpha) \rangle \\
  & = \lim_\alpha \langle \eta \tensor \xi_\alpha, W(\eta \tensor \xi_\alpha) \rangle,
  \end{split}
\]
which means that (\ref{asycoam}) holds.
\par 
(ii) $\Longrightarrow$ (i): Let $a_1, b_1, \ldots, a_\nu, b_\nu \in {\cal C}_0(\G)$, and let $\epsilon > 0$. It is enough to show that there is a vector $\xi \in \Ball(L^2(\G))$ with $\| \xi \| \geq 1 - \epsilon$ such that
\begin{equation} \label{goal2}
  \| M_{a_j,b_j} \circ (W|\xi) - a_j b_j \tensor \xi \|_{{\cal C}_0(\G) \wTensor L^2(\G)_c} < \epsilon \qquad (j=1,\ldots, \nu).
\end{equation}
\par 
Since $\hat{\G}$ is co-amenable, there is a net $( \xi_\alpha )_{\alpha \in \mathbb A}$ of unit vectors in $L^2(\G)$ such that (\ref{asycoam}) holds; it follows easily from (\ref{asycoam}) that
\begin{equation} \label{asycoam2}
  \| \lambda_2(f)\xi_\alpha - \langle f,1 \rangle \xi_\alpha \| \to 0 \qquad ( f\in L^1(\G)).
\end{equation}
For $\alpha \in \mathbb A$, define
\[
  T_\alpha \!: L^1(\G) \to L^2(\G), \quad f \mapsto \lambda_2(f) \xi_\alpha - \langle f, 1 \rangle \xi_\alpha.
\]
The net $( T_\alpha )_\alpha$ lies in $\CB(L^1(\G),L^2(\G)_c)$, is norm bounded, and converges to $0$ pointwise on $L^1(\G)$ by (\ref{asycoam2}). 
\par 
Let $m_0 \in L^1(\G)$ be an arbitrary state, and define, for $j = 1, \ldots, \nu$, matricial subsets $\boldsymbol{K}_j = ( K_{j,n} )_{n=1}^\infty$ of $L^1(\G)$ just as in the proof of Theorem \ref{P1thm}. By Lemma \ref{cunilem} and (\ref{asycoam2}), there is $\alpha_\epsilon \in \mathbb A$ such that
\[
  \sup_{n \in \posints} 
  \sup \left\{ \left\| T_{\alpha_\epsilon}^{(n)} f \right\| : f \in K_{1,n} \cup \cdots \cup K_{\nu,n} \right\} < \frac{\epsilon}{2} 
\]
as well as
\begin{equation} \label{ineq}
  \| \lambda_2(m_0) \xi_{\alpha_\epsilon} - \xi_{\alpha_\epsilon} \| < 
  \frac{1}{\max \{ \| a_1 b_1 \|, \ldots, \| a_\nu b_\nu \| , 1 \}} \frac{\epsilon}{2}.
\end{equation}
Set $\xi := \lambda_2(m_0) \xi_{\alpha_\epsilon}$. It is clear that $\| \xi \| \leq 1$, and by (\ref{ineq}), we also have $\| \xi \| > 1 - \frac{\epsilon}{2} > 1-\epsilon$.
\par 
To prove (\ref{goal2}) holds, first note that
\[
  \begin{split}
  \| M_{a_j,b_j} \circ (W|\xi) - a_j b_j \tensor \xi \| 
  & \leq \| M_{a_j,b_j} \circ (W|\xi) - a_j b_j \tensor \xi_{\alpha_\epsilon} \| 
    + \| a_j b_j \tensor \xi_{\alpha_\epsilon} - a_j b_j \tensor \xi \| \\
  & = \| M_{a_j,b_j} \circ (W|\xi) - a_j b_j \tensor \xi_{\alpha_\epsilon} \| 
    + \| a_j b_j \| \| \xi - \xi_{\alpha_\epsilon} \| \\
  & < \| M_{a_j,b_j} \circ (W|\xi) - a_j b_j \tensor \xi_{\alpha_\epsilon} \| + \frac{\epsilon}{2} 
  \qquad ( j =1, \ldots, \nu),
\end{split}
\]
so that it is sufficient to show that
\begin{equation} \label{subgoal2}
  \| M_{a_j,b_j} \circ (W|\xi) - a_j b_j \tensor \xi_{\alpha_\epsilon} \| < \frac{\epsilon}{2} \qquad (j=1,\ldots, \nu).
\end{equation}
\par 
With $j \in \{1, \ldots, \nu \}$ fixed, observe that
\begin{multline} \label{sups2}
  \| M_{a_j,b_j} \circ (W|\xi) - a_j b_j \tensor \xi_{\alpha_\epsilon} \| \\
  = \sup_{n \in \posints} \sup \{ \| [ (M_{a_j,b_j} \circ (W|\xi)) \eta_{k,l}  
      - \langle \xi_{\alpha_\epsilon}, \eta_{k,l} \rangle a_j b_j ] \|_n : [ \eta_{k,l} ] \in \Ball(\mathbb{M}_n(L^2(\G)_c^\ast)) \}  
\end{multline}
and that the second supremum of the right hand side of (\ref{sups2}) is 
\begin{multline} \label{sup22}
  \sup \{ \| \llangle  [ (M_{a_j,b_j} \circ (W|\xi))\eta_{k,l} 
      - \langle \xi_{\alpha_\epsilon}, \eta_{k,l} \rangle a_j b_j ], [\mu_{\kappa,\lambda}] \rrangle \|_{n^2}: \\
     [ \eta_{k,l} ] \in \Ball(\mathbb{M}_n(L^2(\G)_c^\ast)), \, [ \mu_{\kappa,\lambda} ] \in \Ball(\mathbb{M}_n(M(\G))) \} .
\end{multline}
Then note that, for $\eta \in L^2(\G)^\ast$ and $\mu \in M(\G)$, we have
\[
  \begin{split}
  \lefteqn{\langle (M_{a_j,b_j} \circ (W|\xi)) \eta - \langle \xi_{\alpha_\epsilon}, \eta \rangle a_j b_j, \mu \rangle} \\
  & = \langle \lambda_2(b_j \mu a_j) \xi - \langle b_j \mu a_j, 1 \rangle \xi_{\alpha_\epsilon}, \eta \rangle \\
  & = \langle \lambda_2(b_j \mu a_j \ast m_0) \xi_{\alpha_\epsilon} 
    - \langle b_j \mu a_j \ast m_0, 1 \rangle \xi_{\alpha_\epsilon}, \eta \rangle,
  \end{split}
\]
so that (\ref{sup22}) can also be computed as
\[
  \sup \{ \| [ \lambda_2(b_j \mu_{\kappa,\lambda} a_j \ast m_0) \xi_{\alpha_\epsilon} 
    - \langle b_j \mu_{\kappa,\lambda} a_j \ast m_0, 1 \rangle \xi_{\alpha_\epsilon} ] \|_n : 
  [ \mu_{\kappa,\lambda} ] \in \Ball(\mathbb{M}_n(M(\G))) \}.
\]
We conclude that
\[
  \begin{split}
  \lefteqn{\| M_{a_j,b_j} \circ (W|\xi) - a_j b_j \tensor \xi_{\alpha_\epsilon} \|} & \\ 
  & = \sup_{n \in \posints} \sup \{ \| [ \lambda_2(b_j \mu_{\kappa,\lambda} a_j \ast m_0) \xi_{\alpha_\epsilon} 
    - \langle b_j \mu_{\kappa,\lambda} a_j \ast m_0, 1 \rangle \xi_{\alpha_\epsilon} ] \|_n : 
    [ \mu_{\kappa,\lambda} ] \in \Ball(\mathbb{M}_n(M(\G))) \} \\
  & \leq \sup_{n \in \posints} \sup \{ \| [ \lambda_2(f_{\kappa,\lambda}) \xi_{\alpha_\epsilon} 
    - \langle f_{\kappa,\lambda} , 1 \rangle \xi_{\alpha_\epsilon} ] \|_n : 
    [ f_{\kappa,\lambda} ] \in K_{j,n} \} \\
  & = \sup_{n \in \posints} \sup \left\{ \| [ T_{\alpha_\epsilon} f_{\kappa,\lambda} ] \|_n : 
    [ f_{\kappa,\lambda} ] \in K_{j,n} \right\} \\
  & < \frac{\epsilon}{2}
  \end{split}
\]
for $j =1, \ldots, \nu$, so that (\ref{subgoal2}) holds.
\end{proof}
\begin{remark}
In the proof of (ii) $\Longrightarrow$ (i), we could have chosen the net $( \xi_\alpha )_\alpha$ satisfying (\ref{asycoam}) from $L^2(\hat{\G})_+$. This, however, does not mean that the resulting net satisfying Definition \ref{P2def} belongs to $L^2(\G)_+$. First of all, even though $L^2(\G) = L^2(\hat{\G})$ holds by the definition of $\hat{\G}$, there is no need for $L^2(\G)_+$ and $L^2(\hat{\G})_+$ to coincide (or even be related). Furthermore, even \emph{if} we could pick a net $( \xi_\alpha )_\alpha$ from $L^2(\G)_+$ such that (\ref{asycoam}) holds, then it is not clear that the resulting net for Definition \ref{P2def} would lie in $L^2(\G)_+$ as well.
\end{remark}
\par 
Combining Theorem \ref{P1thm} and \ref{P2thm} and \cite[Theorem 3.2]{BT}, we obtain:
\begin{corollary}
Let $\G$ be a locally compact quantum group with $(P_2)$. Then $\G$ has $(P_1)$.
\end{corollary}
\begin{remarks}
\item We believe that $(P_1)$ and $(P_2)$ are, in fact, equivalent, which---in view of Theorems \ref{P1thm} and \ref{P2thm}---would immediately yield Leptin's theorem for locally compact quantum groups. For a locally compact group $G$, the implication from $(P_1)$ to $(P_2)$ is a straightforward consequence of the elementary inequality
\begin{equation} \label{ineq2}
  \| f - g \|_2^2 \leq \| f^2 - g^2 \|_1 \qquad (f,g \in L^2(G)_+).
\end{equation}
There is a non-commutative variant of (\ref{ineq2}) for von Neumann algebras in standard form (\cite[Theorem IX.1.2(iv)]{Tak2}), however, in order to get from $(P_1)$ to $(P_2)$ in the group case, we have to apply (\ref{ineq2}) to $L^2$-valued, continuous functions on $G$. We thus believe that, in order to derive $(P_2)$ from $(P_1)$ in a general quantum group context, an operator valued version of \cite[Theorem IX.1.2(iv)]{Tak2} is necessary, e.g., in the framework of $\cstar$-valued weights (see \cite{KV0} and \cite[Section 1]{Kus}, for instance).
\item We have not dealt with property $(P_p)$ for locally compact quantum groups for any $p \in [1,\infty)$ other than $1$ or $2$. For any von Neumann algebra $\M$ and $p \in (1,\infty)$, there is a unique so-called non-commutative $L^p$-space $L^p(\M)$ (see \cite{Haa}, \cite{Izu}, and \cite{Ter} for various constructions). For a locally compact quantum group $\G$, one might thus define $L^p(\G)$ as $L^p(L^\infty(\G))$. However, it seems to be unclear, at least for now, how $L^1(\G)$ could be made to act on $L^p(\G)$ in a satisfactory manner that would enable us to even define $(P_p)$ for arbitrary $p$. For a locally compact group $G$, B.\ E.\ Forrest, H.\ H.\ Lee, and E.\ Samei recently equipped $L^p(\hat{G})$ with an $L^1(\hat{G})$-, i.e., $A(G)$-, module structure (\cite{FLS}), and a related, but not entirely identical construction was carried out by the first named author in \cite{Daw}. Both in \cite{FLS} and \cite{Daw}, the non-commutative $L^p$-spaces are obtained through complex interpolation, following \cite{Izu}. It remains to be seen whether the constructions from \cite{FLS} or \cite{Daw} can be extended to general locally compact quantum groups and whether they can be used to define, in a meaningful way, property $(P_p)$ for locally compact quantum groups for arbitrary $p$.
\end{remarks}
\renewcommand{\baselinestretch}{1.0}
\dated
\renewcommand{\baselinestretch}{1.2}
\vfill
\begin{tabbing} 
\textit{Second author's address}: \= Department of Mathematical and Statistical Sciences \kill
\textit{First author's address}:  \> Department of Pure Mathematics \\
                                  \> University of Leeds \\
                                  \> Leeds, LS2 9JT \\
                                  \> United Kingdom \\[\medskipamount]
\textit{E-mail}:                  \> \texttt{matt.daws@cantab.net} \\[\bigskipamount]  
\textit{Second author's address}: \> Department of Mathematical and Statistical Sciences \\
                                  \> University of Alberta \\
                                  \> Edmonton, Alberta \\
                                  \> Canada T6G 2G1 \\[\medskipamount]
\textit{E-mail}:                  \> \texttt{vrunde@ualberta.ca}
\end{tabbing}
\end{document}

%% file: format.tex
\renewcommand{\baselinestretch}{1.2}
\newcommand{\dated}{\mbox{} \hfill {\small [{\tt \today}]}}

%% file: mathdefs.tex
\usepackage{amsmath,amssymb,amsfonts,diagrams}
\newenvironment{keywords}{\noindent\small {\it Keywords\/}:}{\vskip 4pt}
\newenvironment{classification}{\noindent\small 2000 {\it Mathematics Subject
Classification\/}:}{\vskip 12pt}

\newcommand{\comps}{{\mathbb C}}

\newcommand{\posints}{{\mathbb N}}

\newcommand{\tensor}{\otimes}
\newcommand{\ttensor}{\tilde{\otimes}}

\newcommand{\wTensor}{\check{\otimes}}

\newcommand{\cstar}{{C^\ast}}

\newcommand{\id}{{\mathrm{id}}}
\newcommand{\cb}{{\mathrm{cb}}}

\newcommand{\A}{{\mathfrak A}}
\newcommand{\B}{{\mathfrak B}}
\newcommand{\Hilbert}{{\mathfrak H}}
\newcommand{\Kilbert}{{\mathfrak K}}
\newcommand{\M}{{\mathfrak M}}

\newcommand{\CB}{{\cal CB}}

\newcommand{\VN}{\operatorname{VN}}

\newcommand{\varcl}[1]{\overline{#1}}

%% file: thmdefs.tex
\usepackage{amsthm,enumerate}
\theoremstyle{plain}
\newtheorem{theorem}{Theorem}[section]
\newtheorem{lemma}[theorem]{Lemma}
\newtheorem{corollary}[theorem]{Corollary}
\newtheorem{proposition}[theorem]{Proposition}
\theoremstyle{definition}
\newtheorem{definition}[theorem]{Definition}
\theoremstyle{remark}
\newtheorem*{remark}{Remark}
\newtheorem*{example}{Example}
\newtheorem*{rems}{Remarks}
\newtheorem*{exs}{Examples}
\newenvironment{remarks}{\begin{rems}\begin{enumerate}}{\end{enumerate}\end{rems}}

\newenvironment{items}{\begin{enumerate}[\rm (i)]}{\end{enumerate}}
\newenvironment{alphitems}{\begin{enumerate}[\rm (a)]}{\end{enumerate}}